\documentclass[times,preprint]{elsarticle}

\makeatletter
\def\ps@pprintTitle{%
  \let\@oddhead\@empty
  \let\@evenhead\@empty
  \def\@oddfoot{\reset@font\hfil\thepage\hfil}
  \let\@evenfoot\@oddfoot
}
\makeatother

\usepackage{graphicx}
\usepackage{amsthm}
\usepackage{amsmath}
\usepackage{amsfonts}
\usepackage{amssymb}
\usepackage{color}

\usepackage[colorlinks=true,citecolor=blue,urlcolor=blue,linkcolor=blue,pdfpagemode=UseNone]{hyperref}

\theoremstyle{plain}
\newtheorem{theorem}{Theorem}[section]
\theoremstyle{definition}
\newtheorem{definition}[theorem]{Definition}
\theoremstyle{remark}

\newcommand{\cP}{{\mathcal P}}

\begin{document}

\begin{frontmatter}

\title{Rigorous computation of escape times for~parameter~intervals in the quadratic map}

\author[addr1]{Pawe\l{} Pilarczyk\corref{cor1}}
\ead[url]{pawelpilarczyk.com}
\cortext[cor1]{Corresponding author. April 8, 2021}

\author[addr2]{Stefano Luzzatto}
\ead[url]{stefanoluzzatto.net}

\address[addr1]{Faculty of Applied Physics and Mathematics \& Digital Technologies Center,
Gda\'{n}sk University of Technology,
ul.~Gabriela Narutowicza 11/12, 80-233 Gda\'{n}sk, Poland}
\address[addr2]{Abdus Salam
International Centre for Theoretical Physics (ICTP),
Strada Costiera 11, 34151 Trieste, Italy}

\begin{abstract}
We study the quadratic family of one-dimensional maps $f_a (x) = a - x^2$.
We conduct comprehensive numerical analysis of collections of finite orbits of the critical point, computed for \emph{intervals of parameter values} using rigorous numerical methods.
We use the computer to explicitly construct a collection of several thousand parameter intervals, contained in $\Omega=[1.4, 2]$, that are \emph{proved} to have a specific so-called \emph{escape time}, which roughly means that some effectively computed iterate of the critical point taken over all the parameters in that interval has considerable width in the phase space. In particular, we compute a \emph{rigorous lower bound} on this width, in addition to the upper bound.
We investigate the effect of certain constraints imposed on the numerical computations upon the resulting collection of intervals. Additionally, we illustrate and discuss the distribution of the computed intervals in the parameter space.
The purpose of our work is to establish grounds for further numerical computation of a lower bound on the measure of stochastic parameters in $\Omega$.
The source code of the software and the data discussed in the paper are freely available at \href{http://www.pawelpilarczyk.com/quadr/}{http://www.pawelpilarczyk.com/quadr/};
this web page also allows carrying out some limited computations.
The ideas and procedures introduced in the paper can be easily generalised to apply to other parametrised families of dynamical systems.
\end{abstract}


\end{frontmatter}


\section{Introduction}
\label{sec:intro}

The quadratic family of one-dimensional maps

\begin{equation}
\label{eq:quadratic}
f_{a}(x) = a-x^{2}
\end{equation}
is one of the most studied examples of a dynamical system and a landmark of Chaos Theory, due to the extremely rich variety of dynamical behaviours it exhibits. In can be said, however, that notwithstanding the fact that early interest in this family is due to remarkable numerical simulations of the dynamics, most of the deeper results are analytic and probabilistic and do not generally give explicit information about the dynamics of specific parameters in the family. We discuss here, and contribute to, recent work which attempts to bring together sophisticated analytic arguments with rigorous numerical results, in order to obtain more quantitative, and thus more effective, results concerning this family of dynamical systems.

\subsection{Regular and Stochastic Dynamics}

In our investigation of the quadratic family, we restrict the parameters to $a\in \Omega:=[1.4, 2]$, since the dynamics of $f_{a}$ is essentially trivial and well understood for $a\notin\Omega$.
We also restrict the initial conditions $x\in I_{a}$, 
where the interval~$I_{a}$ depends continuously on the parameter $a$ and has the property that $f(I_{a})\subseteq I_{a}$, and that the iterates of all the points $x\notin I_{a}$ converge to $-\infty$. The existence of $I_{a}$ follows by elementary observations and its properties imply that any non-trivial dynamics is contained in $I_{a}$.

There are some very deep abstract results obtained by top mathematicians regarding the quadratic family,
and also extensive numerical studies, carried out mainly by physicists, starting from Feigenbaum.
In spite of these achievements, rigorous numerical results regarding this family are few and far between.
Indeed, numerical analysis of this family of maps is especially challenging due to strong yet non-uniform expansion of the maps
and the properties of chaotic dynamics that are \emph{observed} for many parameters.

We emphasise the word \emph{observed} in the previous sentence in order to indicate that there is a world of difference between phenomenological observation using numerical simulations and mathematically rigorous \emph{proof} of the existence of chaotic dynamics. Indeed, we are especially interested in actually proving the fact that chaotic dynamics is actually present for \emph{many} parameters. The problem is deeper than it may seem. Although the richness of dynamics and the dependence on the parameter are extremely complicated, it is known \cite{AviLyudMel03, Lyu02} that there are only two kinds of dynamics that occur with positive probability in the parameter interval $\Omega$: the dynamics is either \emph{regular}, where $f_{a}$ admits a unique attracting periodic orbit to which Lebesgue almost every $x\in I_{a}$ converges, or the dynamics is \emph{stochastic}, where $f_{a}$ admits a unique invariant probability measure $\mu_{a}$ to which the ergodic averages of Lebesgue almost every point $x\in I_{a}$ converge (in a very ``chaotic'' way, thus the term ``stochastic''). We call the corresponding parameters \emph{regular} and \emph{stochastic} for short.

The difficulty in establishing the set of stochastic parameters is mainly caused by the fact that the set of regular parameters is \emph{open and dense} in $\Omega$ \cite{GraSwi97, Lyu97}, and therefore, spotting a specific stochastic parameter different from $a = 2$ is an extremely rare opportunity. In spite of the fact that the set of stochastic parameters is nowhere dense, it turns out that it has \emph{positive Lebesgue measure} \cite{BenCar85, Jak81}, which shows that it cannot be neglected. Providing an \emph{explicit lower bound} on this measure, however, turns out to be an extremely cumbersome and laborious task. Accomplishing this task is not a matter of simply ``keeping track of constants,'' but it requires reformulation of the original arguments. Successful attempts were made in \cite{Jak01}, where an algorithm was designed towards this purpose, and independently in \cite{LuzTak06}, where an explicit lower bound (though extremely small) was actually obtained. 

An explicit and rigorous lower bound for the set of regular parameters in the closely related \emph{logistic family} \(f_\lambda (x) = \lambda x(1-x) \) was obtained in \cite{TucWil09} through the explicit computation of around 5 million subintervals consisting of regular parameters. The computations required the equivalent of a whole year of CPU time and yielded a set of parameters whose combined measure amounted to little more than 10\% of the length of the entire relevant parameter space.  The logistic family is conjugated to the quadratic family by a smooth map between the corresponding parameter spaces and thus we can obtain a similar estimate for the quadratic family, but this still leaves around 90\% of the parameters completely unaccounted for, with no rigorous nor heuristic argument currently able to establish whether most of these parameters are regular or stochastic.

\subsection{Escaping Intervals}

Notwithstanding the complexity of the task of bounding the measure of stochastic parameters from below due to its topological nowhere dense structure, it is arguably an ultimately more effective way to approach the problem than trying to explicitly compute infinitesimally small intervals of regular parameters. Indeed, a very promising strategy in this direction is to take advantage of the powerful analytic arguments which have been developed to prove the positive probability of the set of stochastic parameters, such as in the pioneering papers \cite{BenCar85, Jak81} and later generalisations \cite{LuzTuc99, LuzTak06, LuzVia00, PRV}, and to develop computer-assisted methods that would allow effective computation of rigorous estimates for certain constants that can be later plugged in into analytical arguments. More specifically, the setup in \cite{LuzTak06} does not use any computer-assisted calculations to get a first ever explicit lower bound on the measure of the set of stochastic parameters, but reformulates the analytic arguments in such a way as to clarify the \emph{explicitly computable} quantities which are required as inputs to the analytical arguments in order to obtain rigorous and explicit bounds; in particular, its setting takes advantage of the special characteristics of the parameter value \( a = 2 \) to compute the required quantities analytically. As a consequence, the approach introduced in \cite{LuzTak06} provides a possible roadmap towards full implementation of such a strategy.

In particular, one of the most important quantities that must be found is an \emph{escape time} defined below (see also \cite{BenCar91}). Let $\omega\subseteq \Omega$ be an arbitrary parameter interval and let $c$ denote the critical point $0$ of $f_a$; notice that the critical point is the same for all parameter values. For each $n\geq 0$, we let 

\begin{equation}\label{eq:omegan}
c_{n}(a):= f_{a}^{n}(f_{a}(c)) \qquad \text{ and } \qquad 
\omega_{n}:=\{c_{n}(a): a\in \omega\}.
\end{equation}
Note that the critical value $c_{0}(a)$ equals $a$; therefore, $\omega_{0}$ coincides numerically with $\omega$ although it lives in phase space rather than parameter space. For $n\geq 1$, $c_{n}(a)$ is simply the $n$-th image of the critical value, and $\omega_{n}$ is the interval given by the $n$-th images of the critical values for all the parameters $a\in \omega$. 

We make a choice of a constant $\delta>0$ that defines the \emph{critical neighbourhood}

\begin{equation}
\label{eq:Delta}
\Delta:=(-\delta, \delta). 
\end{equation}

\begin{definition}[escape time]
\label{def:escape}
$N$ is called an \emph{escape time} for $\omega$ if the following holds:
\begin{equation}
\label{eq:escape}
\omega_{i}\cap \Delta=\emptyset \quad \text{for all $i \in \{0, \ldots, N-1\}$}, \quad \text{and}  \quad |\omega_{N}|\geq \sqrt{\delta}. 
\end{equation}
\end{definition} 
The existence of an escape time for $\omega$ roughly means that some iterate of the critical point taken over all the parameters in $\omega$ has considerable width in the phase space, and thus can serve as a starting point for proving that it contains a positive measure of stochastic parameters. Indeed, the existence of an escape time for a sufficiently large value of \( N \)  is  arguably the key ingredient in the proof of positive measure of stochastic parameters. It therefore becomes an important problem to be able to explicitly compute intervals of parameters which have escape times and to compute the corresponding escape times.

\subsection{Algorithms for the Computation of Escaping Intervals}

In our previous paper \cite{GolKouLuzPil20}, we developed an algorithm that allowed us to effectively construct almost $1.5$ million parameter intervals $\{ \omega^j \}$ for which we proved the existence of a number $N(\omega^j) \geq N_0 := 25$ that was an escape time for the corresponding $\omega^j$ with $\delta = 10^{-3}$. We found out that these intervals covered almost $90\%$ of $\Omega$, and we conducted some basic analysis of selected features of the constructed intervals.

Let us briefly recall the algorithm; the reader is referred to \cite{GolKouLuzPil20} for the details, and to \cite{WT2011} for general introduction to rigorous numerics. We first fix some radius $\delta > 0$ for the critical neighbourhood $\Delta$. Then we split $\Omega$ uniformly into $u \geq 1$ subintervals and put them all into a queue $Q$. We fix a desired number $N_0 > 0$ of iterates in computing subsequent $\omega_i$ for $i = 1, 2, \ldots$, and an upper bound $N_{\max} \geq N_0$ on the number of the iterates to consider. The justification for setting $N_{\max}$ is the following. If an interval $\omega \subset \Omega$ can be iterated many times without hitting $\Delta$ then it is likely that the intervals $\omega_i$ have been trapped in a neighbourhood of an attracting periodic orbit, so there is no hope to get a reasonable escape time for it, because the intervals $\omega_i$ shrink in the subsequent iterations.

In the main loop of the algorithm, we repeatedly take an interval $\omega$ from the queue $Q$, and compute consecutive iterates $\omega_i$, as defined in \eqref{eq:omegan}, for $i = 1, 2, \ldots$, until one of the following cases labelled (P1)--(P3) happens: either (P1) $\omega_i \cap \Delta \neq \emptyset$, in which case we say that \emph{$\omega_i$ hits $\Delta$} (or we say that \emph{$\omega$ hits $\Delta$ after $i$ iterations}), or (P2) a numerical problem appears, which is most commonly the inability to prove that $c_n'$ is of constant sign in $\omega$, an assumption necessary to continue the iterations, or, eventually, (P3) $i > N_{\max}$. In the case (P1), if $i \geq N_0$ and $|\omega_i| \geq \sqrt{\delta}$ then we add $\omega$ to the collection of ``successful'' intervals $\cP^+$ for which we have found a large enough escape time. In the other cases, we split the interval $\omega$ into smaller parts and put (some of) them back into the queue for another attempt of iterating, or we put them in the set $\cP^-$ of ``failure'' items. A reason for putting $\omega$ or its part in $\cP^-$ is either no hope for successful completion of (P1) with $i \geq N_0$, or the size of $\omega$ below certain fraction $w > 0$ of the width of $\Omega$, or the growth of overestimates exceeding the precision of numerics (taken as $p = 250$ bits per number in \cite{GolKouLuzPil20}). In particular, in the case (P1), if $i < N_0$ or $|\omega_i| < \sqrt{\delta}$, the interval $\omega$ is split into up to three parts in order to exclude the portion of $\omega$ whose $i$-th iterate intersects $\Delta$; this portion is estimated using the bisection method with the number $s > 0$ of steps (set to $s = 40$ in~\cite{GolKouLuzPil20}).

Choosing conservative values for $N_{\max}$ and $w$ ensures that the algorithm completes its work within reasonable time; for example, in \cite{GolKouLuzPil20} we chose $N_{\max} = 200$ and $w = 10^{-10}$, which makes the program complete the computations for $N_0 = 25$ in about 1 hour 10 minutes, using up to 767~MB of RAM on our reference system described at the end of Section~\ref{sec:intro}.

Alternatively, one can consider some other stopping conditions, like the maximum number $i_{\max}$ of intervals processed by the algorithm, or the maximum size of the queue $Q$ (to avoid memory overflow). It is important to note that, together with each interval $\omega$ stored in the queue $Q$, we also remember the number of times the interval was successfully iterated before it was put in the queue. Then we may fix some number $N_{\min} > 0$ and stop pulling intervals from the queue when all the intervals that are left there have been iterated successfully at least $N_{\min}$ times.

\subsection{Statistics of Escaping Intervals}

As described above, the results contained in \cite{GolKouLuzPil20} depend on the choices made for several quantities that play a role in the construction. It is natural to ask how sensitive or robust the results are with respect to these choices. The purpose of the current paper is to investigate the effects on the results caused by making specific choices of the constants and other constraints in the algorithms, and also to get insight into the properties of the collections of parameter intervals constructed by the algorithm. For that purpose, we conducted several rounds of computations with various settings, and we report on our findings in the next sections. Specifically, in Section~\ref{sec:first}, we conduct the analysis of the numbers of iterates and the sizes of the images of subintervals of $\Omega$ at the moment of the first time their iterations hit $\Delta$. In Section~\ref{sec:chopping}, we provide insight into the contribution of intervals of various sizes to the total measure of constructed intervals with an escape time bounded from below by certain values of $N_0$. Finally, in Section~\ref{sec:bisection}, we address the delicate numerical question of choosing the number of bisection steps to be made while chopping an interval whose iterate hits $\Delta$.

In a nutshell, our investigation shows that while the spcific escape intervals and times which arise from the construction may indeed be sensitive to the choice of constants, the overall statistics are instead quite stable. Together with the fact that we are able to include most parameter intervals in escaping intervals this indicates that the method is quite effective. We emphasise that we are not just describing results of numerical \emph{simulations}. These are all high-precision calculations conducted with controlled rounding directions (interval arithmetic), rigorous derivative estimates and monotonicity arguments, and thus constituting collections of authentic facts actually \emph{proved} using the rigorous computational methods.

All the computations described in the paper were conducted using a $64$-bit Debian Linux 10.8 server with the Intel\textsuperscript{\tiny\textregistered} Core\textsuperscript{\tiny\texttrademark}2 Duo Processor E8400 running at 3~GHz, and all the computing times given in the paper refer to this system.

The source code of the software and the data discussed in the paper are freely available at the website \cite{software} that additionally allows carrying out some limited computations.


\section{Parameter intervals at first encounter with the critical neighbourhood $\Delta$}
\label{sec:first}

Consider a small interval $\omega \subset \Omega$ and its consecutive iterates $\omega_1, \omega_2, \ldots$ The interval $\omega$ can sometimes be iterated for a long time before it hits the critical neighbourhood $\Delta$, or this can happen after just a few iterations. Its size might have grown considerably by that time, or might be very small.

In this section, we analyse selected features of collections of small intervals of parameters at their first encounter with $\Delta$ during the iteration process. In order to collect data for the analysis, we split $\Omega$ uniformly into a large number of small intervals, and we compute $\omega_1, \omega_2, \ldots$ for up to $N_{\max} := 100$ iterations, until the first number $N(\omega)$ is found for which $\omega_N \cap \Delta \neq \emptyset$. We fix the precision of the numerics at $200$ bits, which yields the accuracy at the order of $10^{-60}$. We try $\delta \in \{ 10^{-2}, 10^{-3}, \ldots, 10^{-7} \}$, and we split $\Omega$ uniformly into $u \in \{ 60, 600, \ldots, 600000 \}$ subintervals in order to get a wide range of data for comparison; since $|\Omega| = 0.6$, these splittings result in subintervals of size $10^{-2}, \ldots, 10^{-6}$.

\begin{figure}[htbp]
\centerline{\includegraphics[width=10cm]{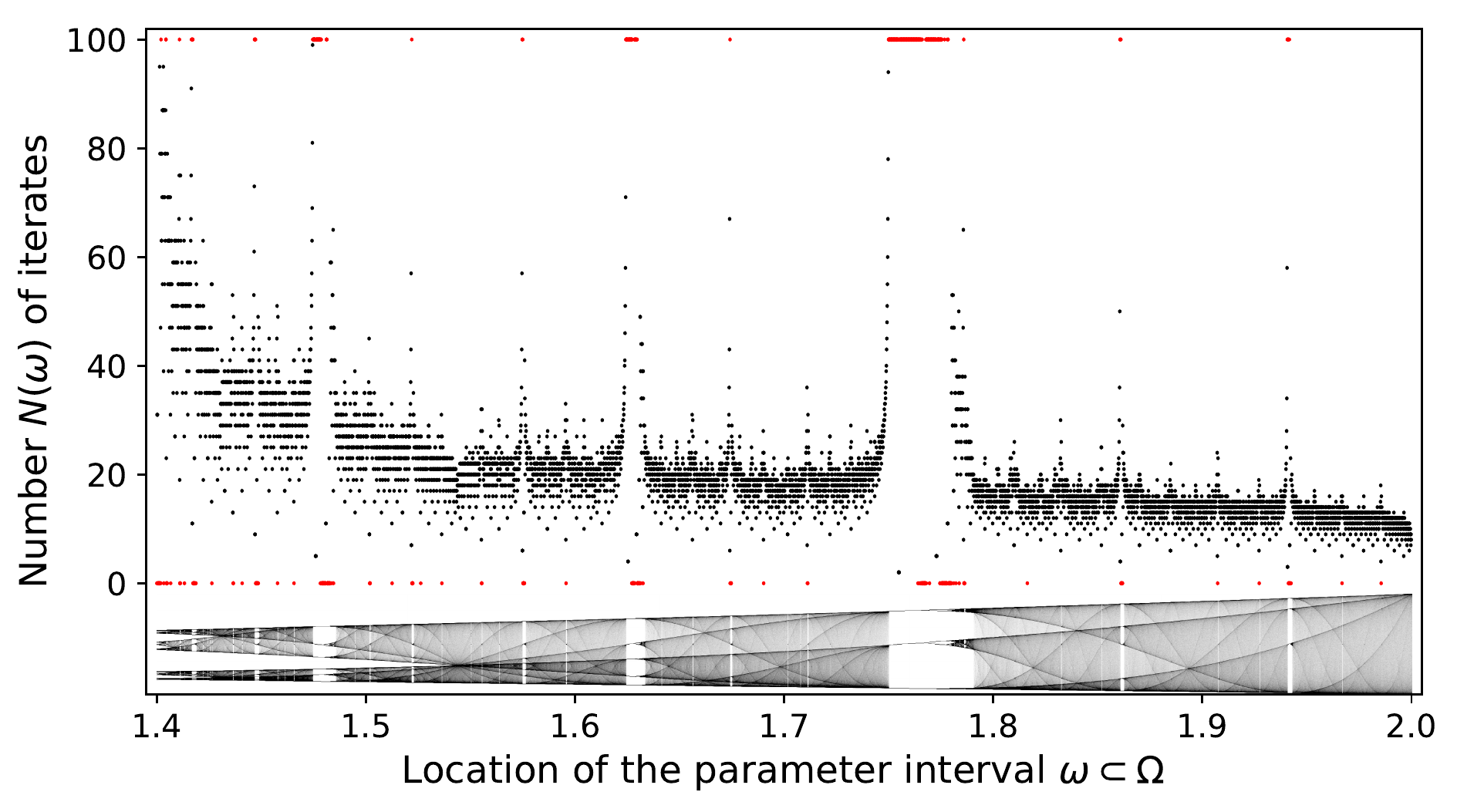}}
\caption{\label{fig:noChNum43} The location and the number $N(\omega)$ of iterates of $u = 6{,}000$ intervals $\omega \subset \Omega$ in a uniform subdivision of $\Omega$ at which the first encounter with $\Delta = (-\delta, \delta)$ was observed for $\delta = 10^{-3}$. The $243$ intervals for which a numerical problem occurred (inability to prove the constant sign of the derivative of $c_i$) are marked in red at level $0$. The $324$ intervals that could be iterated $100$ or more times are marked in red at level $100$. The bifurcation diagram of the quadratic map is drawn along the parameter axis for reference.}
\end{figure}

Figure~\ref{fig:noChNum43} shows the number of iterates until hitting $\Delta$ for intervals in a uniform subdivision of $\Omega$. In the area of high expansion (especially around $a = 2$), the number of iterates is typically smaller than elsewhere. Larger numbers of iterates are achieved close to attracting periodic orbits; however, such intervals are typically useless for a high escape time, because they shrink while following the stable orbit, and thus never achieve the desired width. The shape of the graph is in clear correspondence with the bifurcation diagram.

\begin{figure}[htbp]
\centerline{\includegraphics[width=10cm]{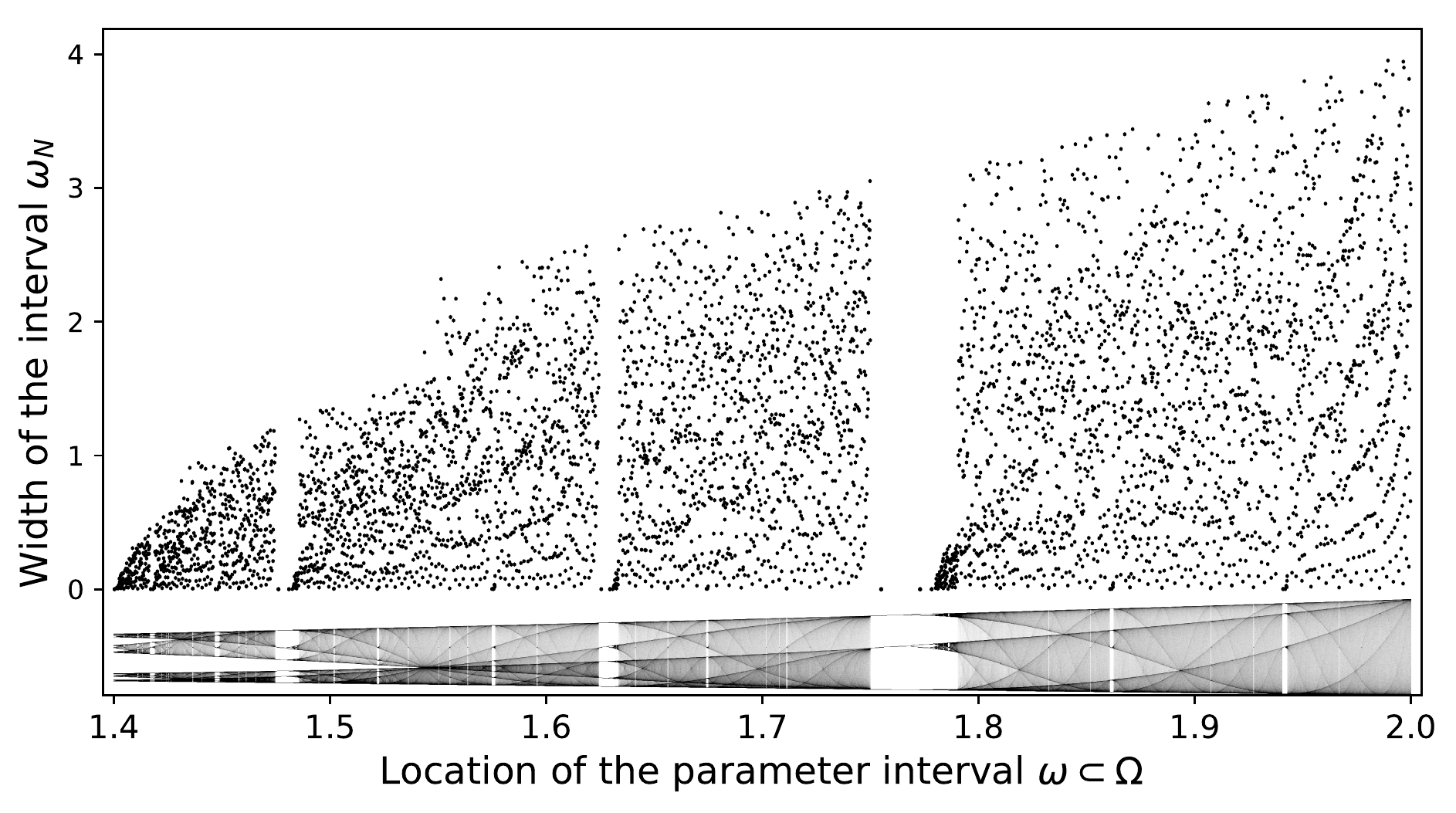}}
\caption{\label{fig:noChWidth43} The widths of the intervals $\omega_{N}$ at the first iterate $N(\omega)$ at which they intersect $\Delta$, computed for $6{,}000$ intervals in a uniform subdivision of $\Omega$. The bifurcation diagram of the quadratic map is drawn along the parameter axis for reference.}
\end{figure}

\begin{figure}[htbp]
\centerline{\includegraphics[width=10cm]{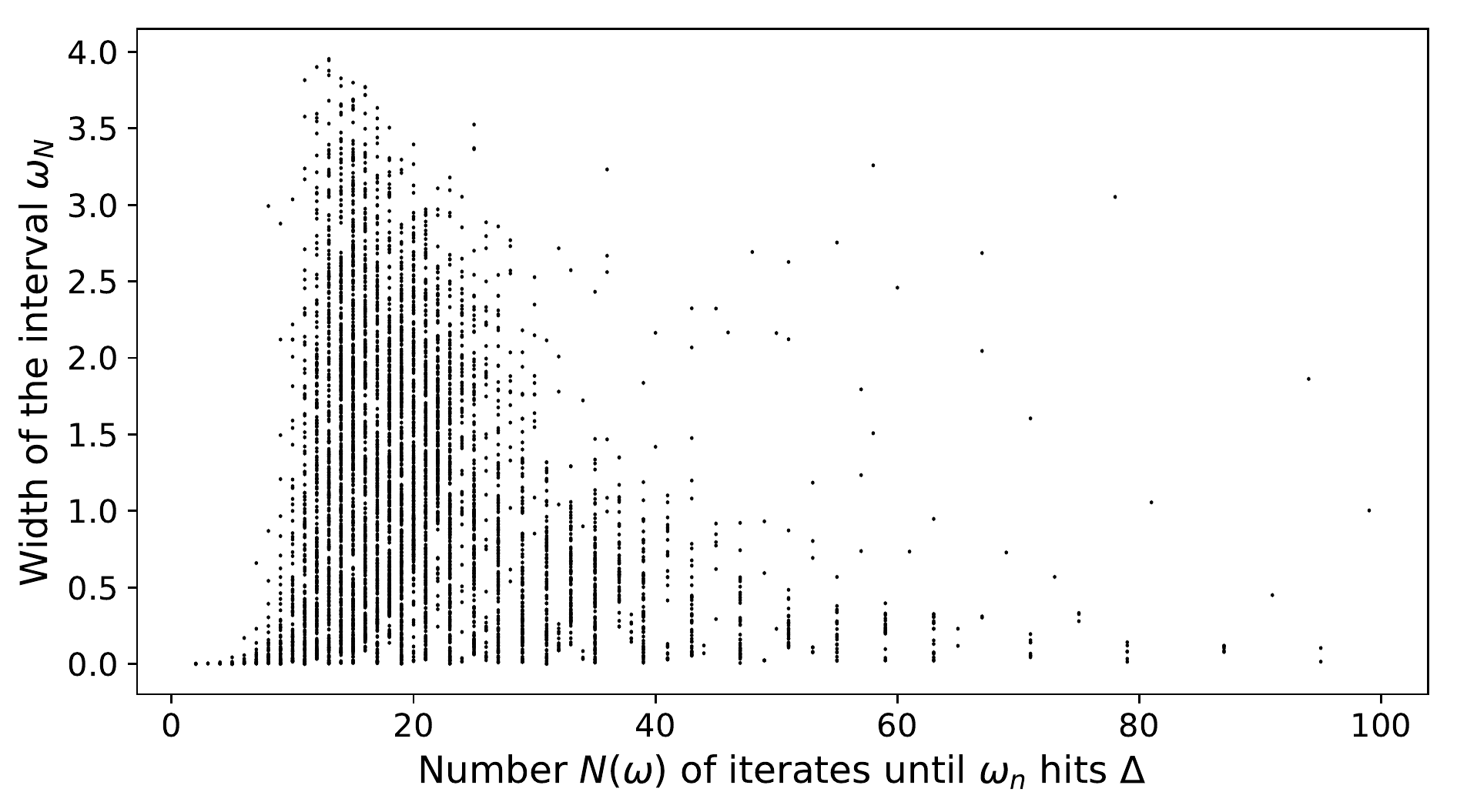}}
\caption{\label{fig:noChPairs43} The widths of the intervals $\omega_{N}$ against the number $N(\omega)$ of iterates after which they hit $\Delta$, computed for $6{,}000$ intervals in a uniform subdivision of $\Omega$.}
\end{figure}

Figure~\ref{fig:noChWidth43} illustrates the widths of the images of the same intervals of parameters as considered previously. It is interesting to see that there is no apparent relation between the location of the interval and the width attained when hitting $\Delta$. In fact, the figure suggests that essentially all possible widths within some range are obtained in all the parameter regions, except for the periodic windows. Indeed, analogous graphs obtained for higher numbers $u$ (not inlcuded in the paper) turn out to be very similar in shape, except these areas are ``darker,'' almost completely filled with data points.

A relation between the number of iterates $N(\omega)$ and the width of $\omega_N$ is illustrated in Figure~\ref{fig:noChPairs43}. An important observation is that, at this size of the parameter intervals, the widest images at time of intersecting $\Delta$ appear at iterations around $12$--$14$, then a gradual decrease in the highest value can be observed, and there are virtually no wide images for high numbers $N(\omega)$. Indeed, this seems to be caused by the fact that strong expansion of the map does not allow iterating the intervals for a long time, except when the iterates are wandering close to a stable periodic orbit, in which case, obviously, the size of the interval shrinks.

\begin{figure}[htbp]
\centerline{\includegraphics[width=10cm]{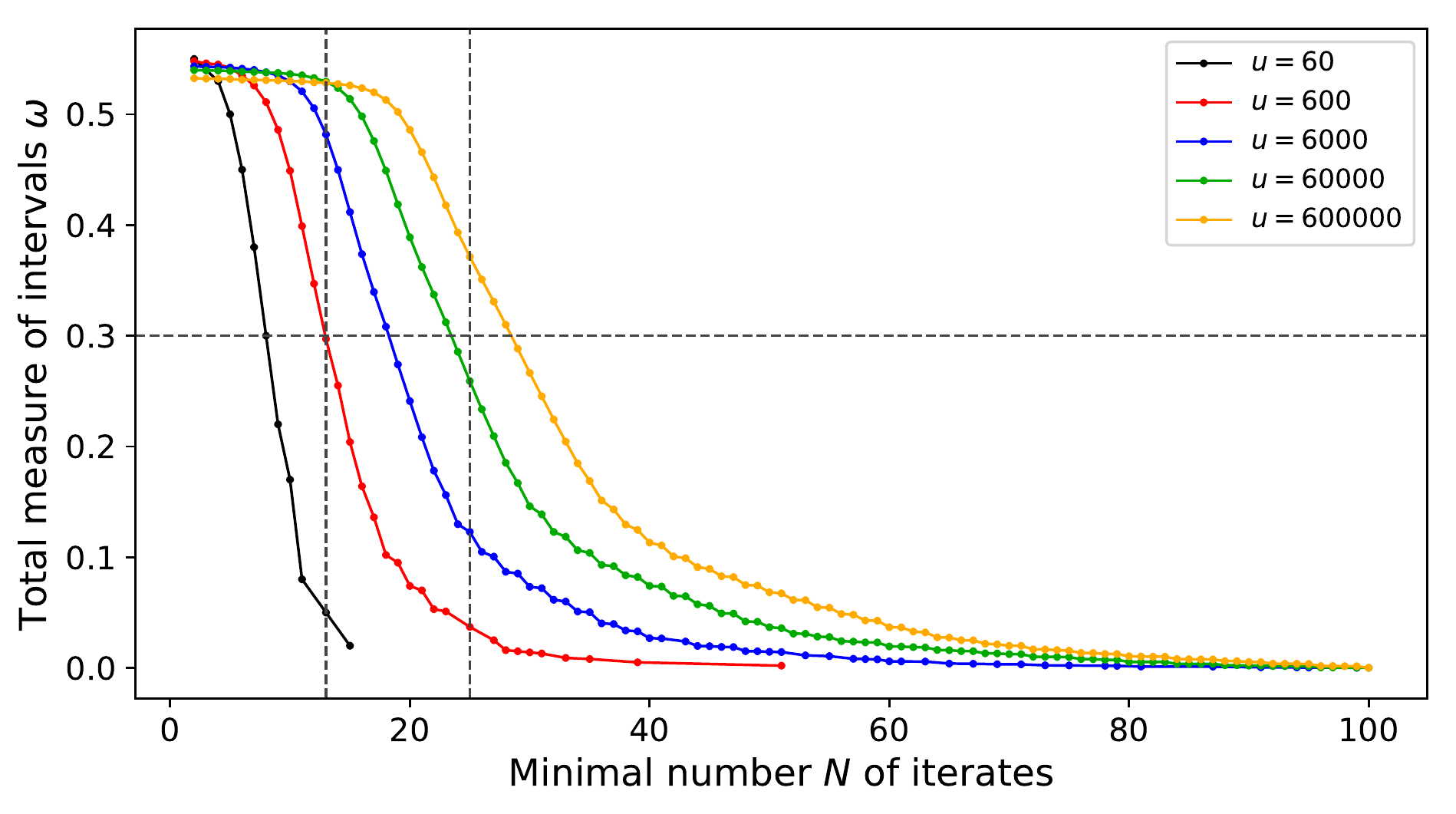}}
\caption{\label{fig:noChNumMin_d3} The total measure of the intervals $\omega$ that can be iterated at least $N$ times until they hit $\Delta = (-\delta,\delta)$, with $\delta = 10^{-3}$, plotted against the number $N$. Five curves are shown, computed for a uniform subdivision of $\Omega$ into five different numbers $u = 60, \ldots, 600000$ of intervals. The values that can be read out from the graph at the grey dashed lines are discussed in the text.}
\end{figure}

We are in fact most interested in the total measure of intervals with a sufficient escape time. Figure~\ref{fig:noChNumMin_d3} shows the total measure of intervals that can be iterated at least a certain number of times, depending on the initial width of the intervals, which corresponds to the number $u$ of the intervals in the uniform subdivision of $\Omega$. For example, the number $u = 600$ implies subdividing the interval $\Omega = [1.4, 2]$ into $600$ subintervals of size approximately $10^{-3}$.

If one looks at the first vertical dashed line in Figure~\ref{fig:noChNumMin_d3} for $N=13$, one can learn from the intersection points of this line with the five curves that the measure of subintervals $\omega \subset \Omega$ that hit $\Delta$ after at least $13$ iterations is about $0.05$ if we consider intervals of size $|\omega| = 10^{-2}$ (a point on the black curve), or almost $0.3$ if we take intervals of length $|\omega| = 10^{-3}$ (a point on the red curve), or about $0.48$ if we take one more order in resolution to $|\omega| = 10^{-4}$ (a point on the blue curve), and almost $0.53$ if $|\omega| = 10^{-5}$ or $|\omega| = 10^{-6}$ (the green and yellow curves). The situation is much less optimistic for the second vertical line, drawn for $N=25$, where there are no intervals of size $|\omega| = 10^{-2}$ that can be iterated that many times (the black line drops to $0$ well before $N=25$), less than a half of $\Omega$ is covered with intervals as small as $10^{-5}$, and the measure of the smallest intervals considered here (of size $10^{-6}$) is only about $0.371$. We are going to compare these results with the ones shown in Figure~\ref{fig:measureW}, where we chop intervals that hit $\Delta$ and continue iterating further.

On the other hand, if one looks at the horizontal dashed line in Figure~\ref{fig:noChNumMin_d3}, one can see that the measure of $0.3$ (which is $50\%$ of the measure $0.6$ of the whole of $\Omega$) comprises of intervals of size $|\omega| = 10^{-2}$ that can be iterated $8$ or more times, while almost the same measure is achieved by intervals of size $|\omega| = 10^{-3}$ that can be iterated at least $13$ times, and eventually by intervals of size $|\omega| = 10^{-6}$ that hit $\Delta$ after at least $28$ times.

The first observation is that the finer the subdivision of $\Omega$, the higher the number of times the small intervals can be iterated. Another interesting feature is the linear drop in the measure as a function of the number $N$ of iterates, observed in the middle part of the graphs. This means, for example, that if we require $3$ more iterates then the drop in the measure would be approximately $3$ times more severe than the decrease caused by requiring only $1$ extra iterate. A less optimistic observation, however, is the constant distance between the curves in the horizontal direction. Each next curve is obtained by splitting $\Omega$ into $10$ times more subintervals. Analysing this distance shows that, while preserving the obtained measure, we must make exponentially higher computational effort in order to get a few more iterates completed.

Although we do not illustrate it here in a separate figure, we would like to mention the fact that the graphs obtained for $\delta$ set to the different values in $\{10^{-2}, 10^{-3}, \ldots, 10^{-7}\}$ are extremely similar to each other, with some small but noticeable differences found at $\delta = 10^{-2}$ only, especially for low values of $N$. An important conclusion from this observation is that the intuition suggesting that decreasing $\delta$ might help the iterated intervals ``avoid'' the smaller critical neighbourhood $\Delta$, is in fact misleading. This is an important point, because decreasing $\delta$ is known to detrimentally affect certain other important estimates (see \cite{DayKokLuzMisOkaPil08,GolLuzPil16}) while possibly help others. Our observation indicates that the actual size of $\Delta$ does not in fact considerably affect the process of iterating the intervals.

\begin{figure}[htbp]
\centerline{\includegraphics[width=10cm]{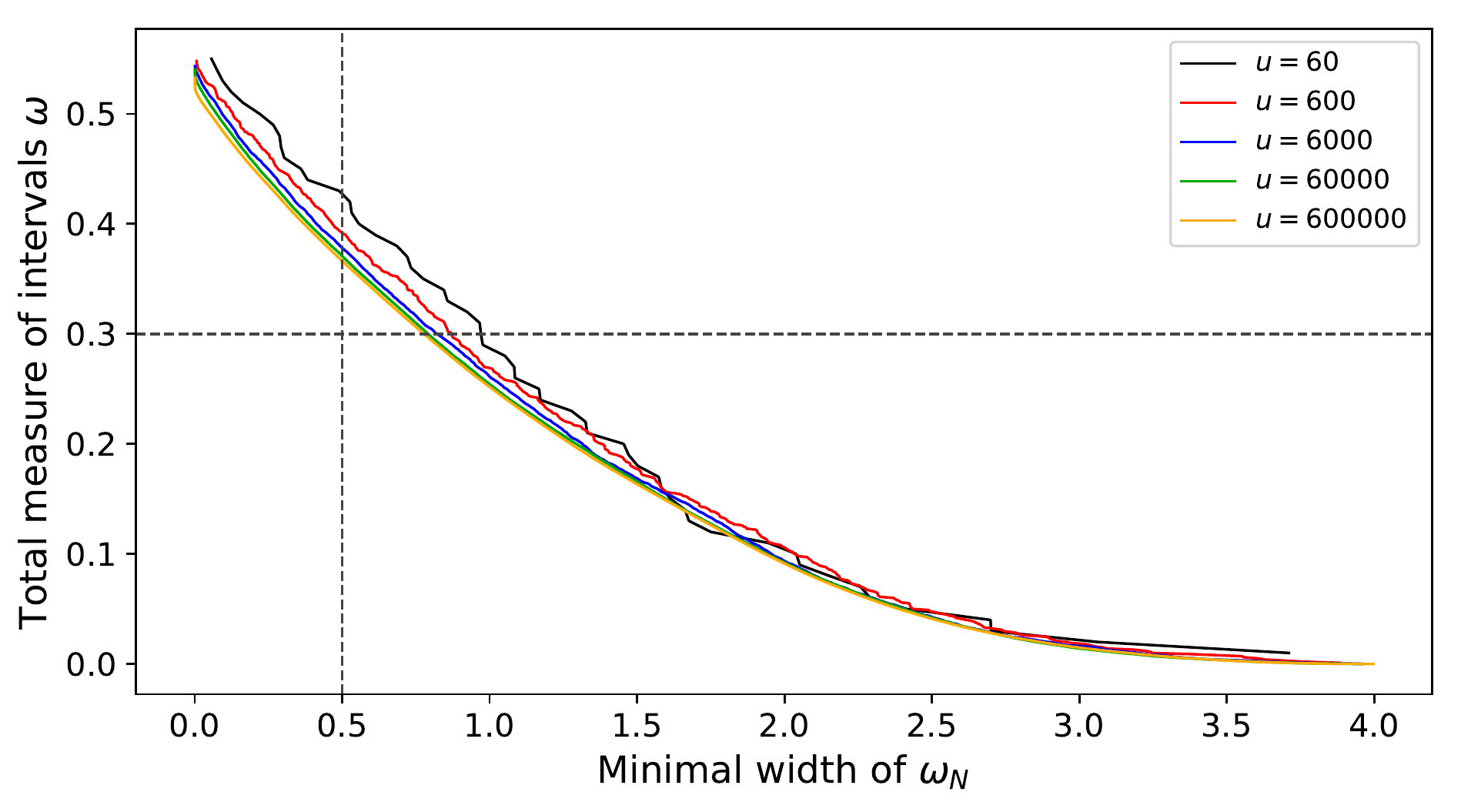}}
\caption{\label{fig:noChWidthMin_d3} The total measure of the intervals $\omega$ whose iterates are of at least the given length at time of hitting $\Delta$, computed with $\delta = 10^{-3}$, computed for the uniform subdivision of $\Omega$ into a few different numbers $u = 60, \ldots, 600000$ of intervals.}
\end{figure}

Figure~\ref{fig:noChWidthMin_d3} shows another feature of intervals with an escape time: the distribution of widths of their images at the time of hitting $\Delta$. If one looks at the vertical dashed line, one can learn from this graph that the measure of subintervals $\omega \subset \Omega$ whose width at time of hitting $\Delta$ is at least $0.5$ (which is indeed a macroscopic size) is between about $0.36$ and $0.43$, depending on $u$. On the other hand, if one looks at the horizontal dashed line, one can see that the measure of $0.3$ (which is $50\%$ of the measure $0.6$ of the whole of $\Omega$) comprises of intervals that grow to the size of $0.77$ or more, or even to $0.97$ or more, depending on~$u$. An observation that evokes mixed feelings is that the widths are larger for smaller values of $u$, that is, for a coarser subdivision of $\Omega$. Unfortunately, this does not indicate that starting with larger subintervals $\omega \subset \Omega$ is beneficial, because the number of iterates until such intervals hit $\Delta$ is then considerably lower, which undermines their usefulness for the escape time condition with large $N_0$. However, the differences in the total measure for different values of $u$ are not as profound as the differences in the number of iterates $N(\omega)$.

An overall conclusion that can be drawn from the analysis conducted in this section is the following. Since the sizes of images of the intervals $\omega \subset \Omega$ are macroscopic in comparison to $\Delta$, shifting $\omega$ a little, shrinking it slightly, or chopping into smaller \mbox{intervals} might make it possible to iterate the interval further, until it intersects $\Delta$ again. With confirmed growth of the image, this would give rise to an escape time with a higher value of $N_0$. Therefore, instead of stopping the iterations at the first encounter with $\Delta$ and discarding intervals for which the number of iterates is insufficient, a considerably more productive approach is to define the subintervals of $\Omega$ dynamically, by means of iterating initially chosen subintervals and chopping them to remove the portion that falls into $\Delta$. In this way, only a small portion of the interval $\omega \subset \Omega$ is chopped off each time, and the majority of the interval is retained for further processing. We discuss the results of such an approach in Section~\ref{sec:chopping}, and we show that they are indeed more appealing.


\section{Chopping parameter intervals and iterating further}
\label{sec:chopping}

In \cite[Fig. 2]{GolKouLuzPil20}, we showed the numbers of intervals of different sizes and the percentage of the measure of $\Omega$ they occupy, obtained after having completed the algorithm described in Section~\ref{sec:intro} with $N_0 = 25$, $\delta = 10^{-3}$, $w = 10^{-10}$, and $p = 250$. In particular, we found out that over $50\%$ of the measure of parameter intervals that have an escape time of at least $N_0$ (see Definition~\ref{def:escape}) was accounted for by a small number (about $7{,}000$) of relatively large intervals (of size above $10^{-5}$), while small intervals of size below $10^{-7}$ contributed less than $6\%$ to the measure, even though their number was huge (above $1.1$ million).

\begin{figure}[htbp]
\noindent
\includegraphics[scale=.32]{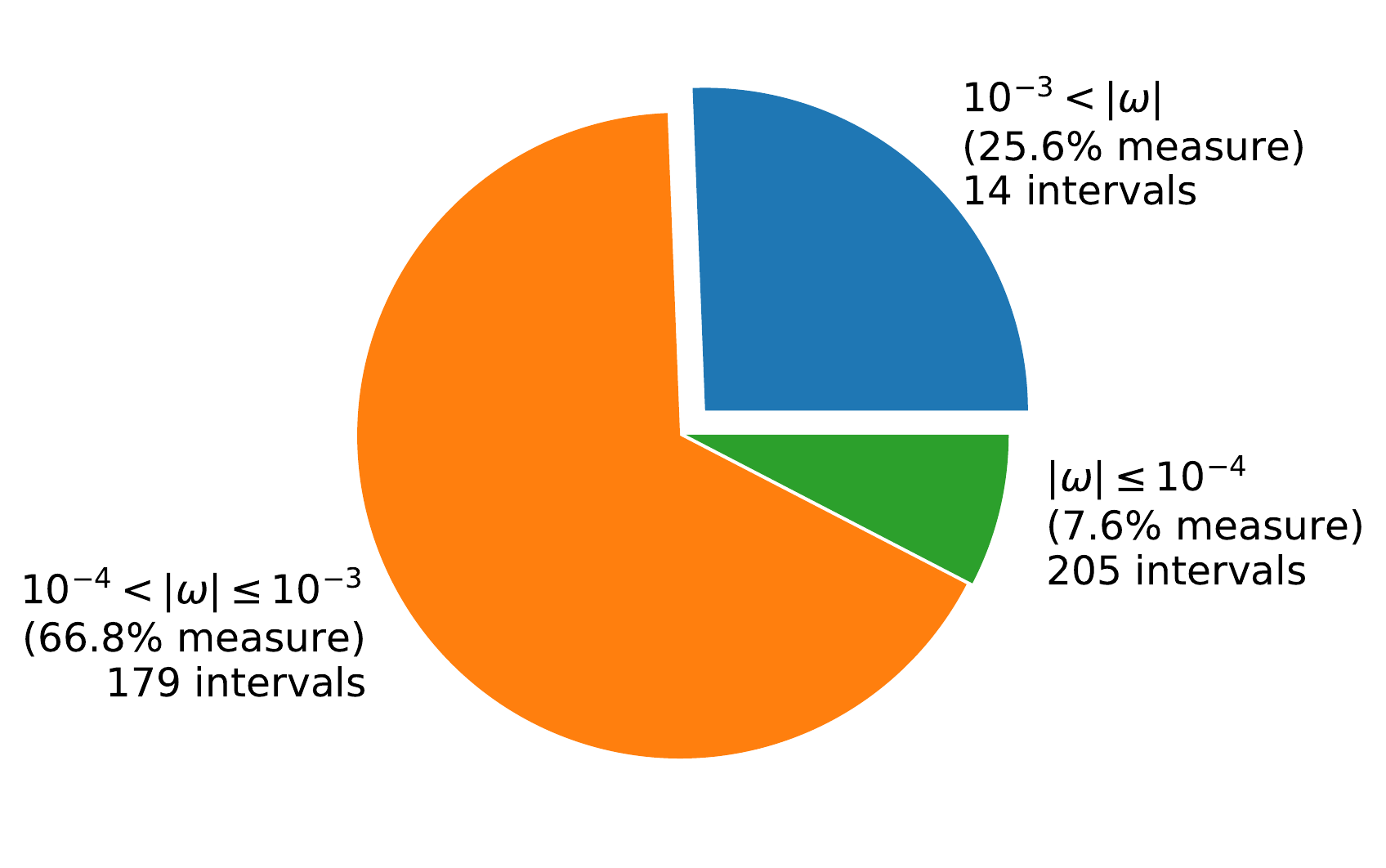}
\includegraphics[scale=.32]{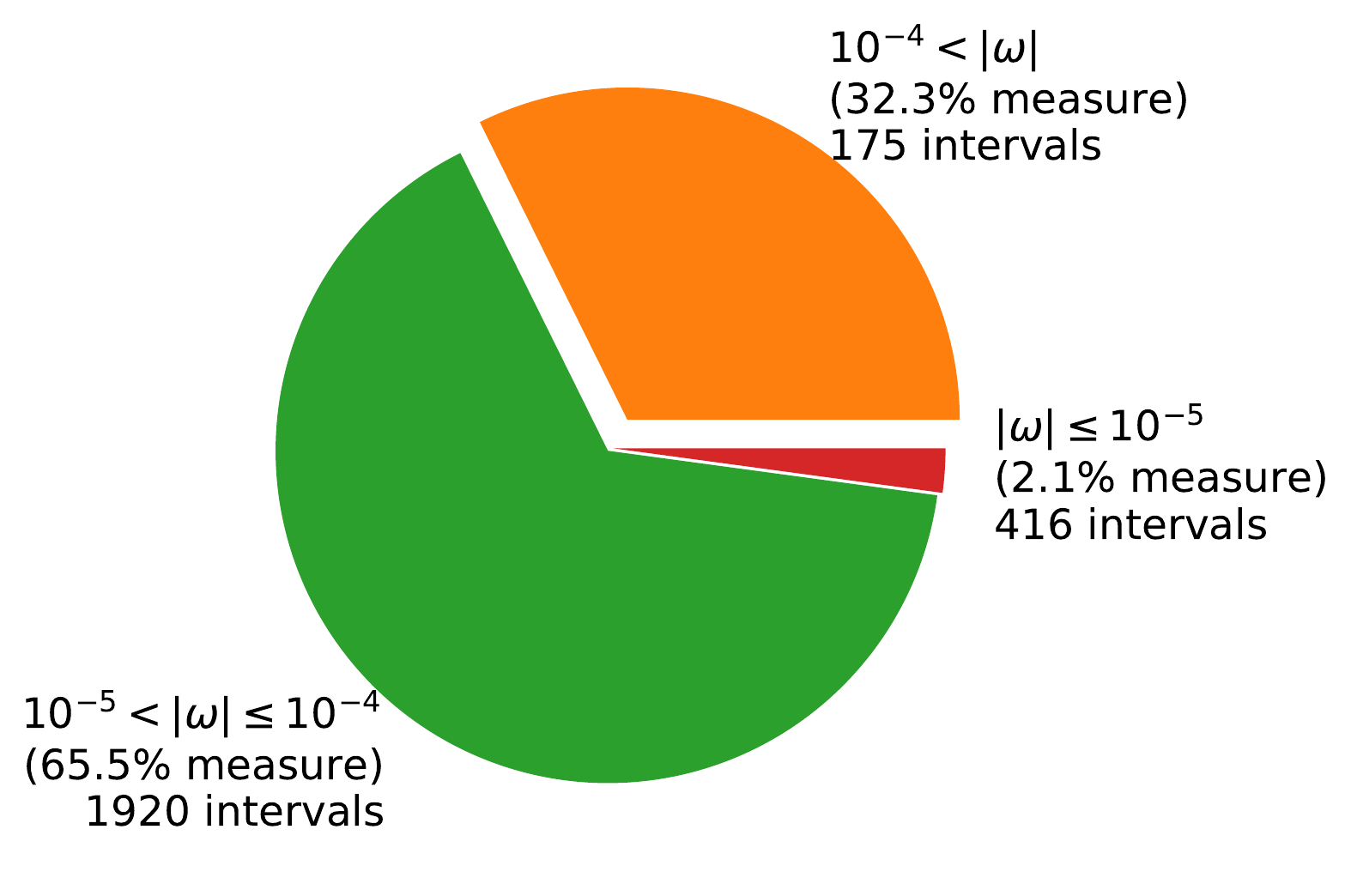}

\centerline{\hfill $\Omega_1 = [1.4,1.5]$, $\mu_1 \approx 85.9\%$ \hfill
\hfill $\Omega_2 = [1.5,1.6]$, $\mu_2 \approx 96.5\%$ \hfill}

\vskip 12pt

\noindent
\includegraphics[scale=.32]{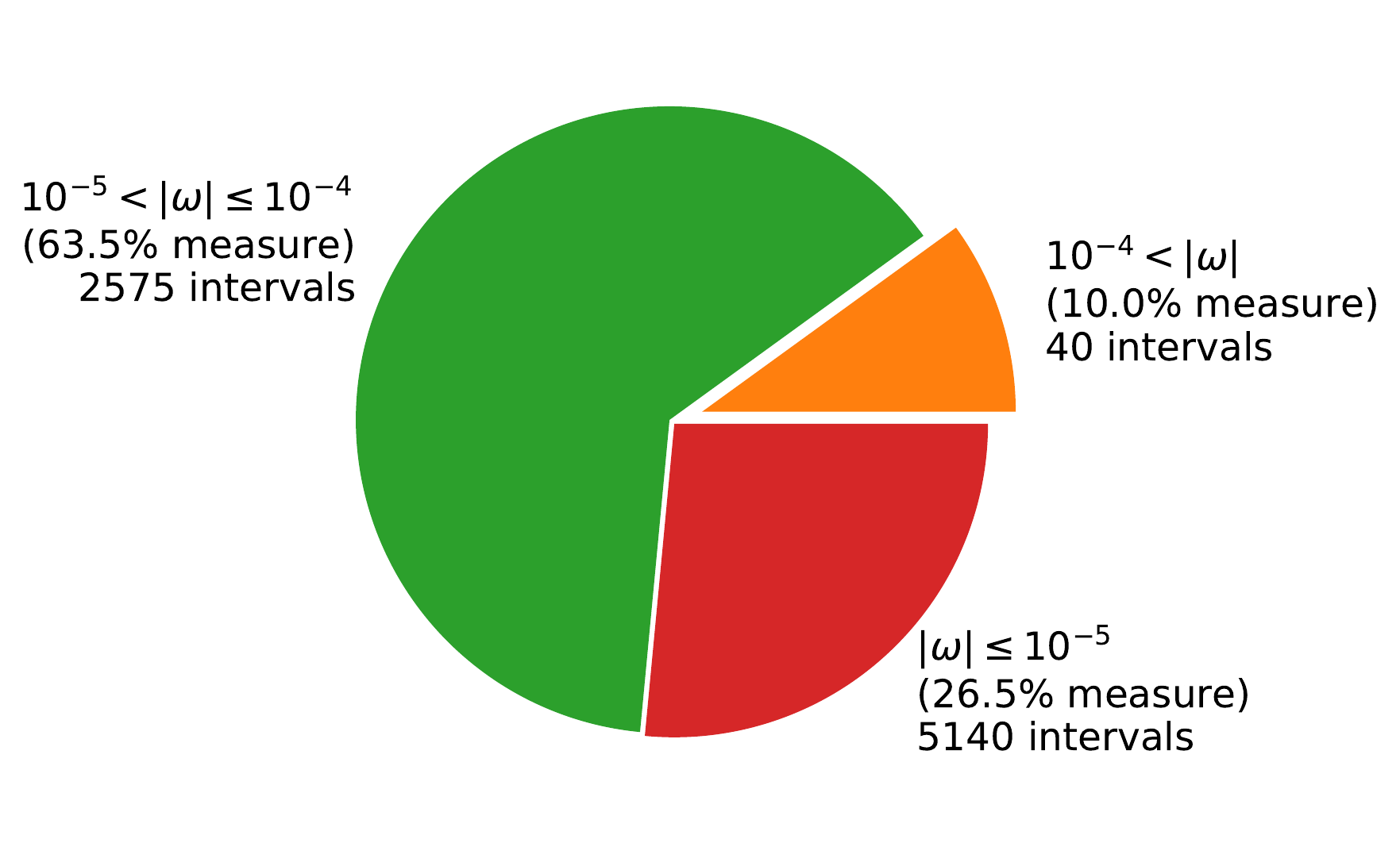}
\includegraphics[scale=.32]{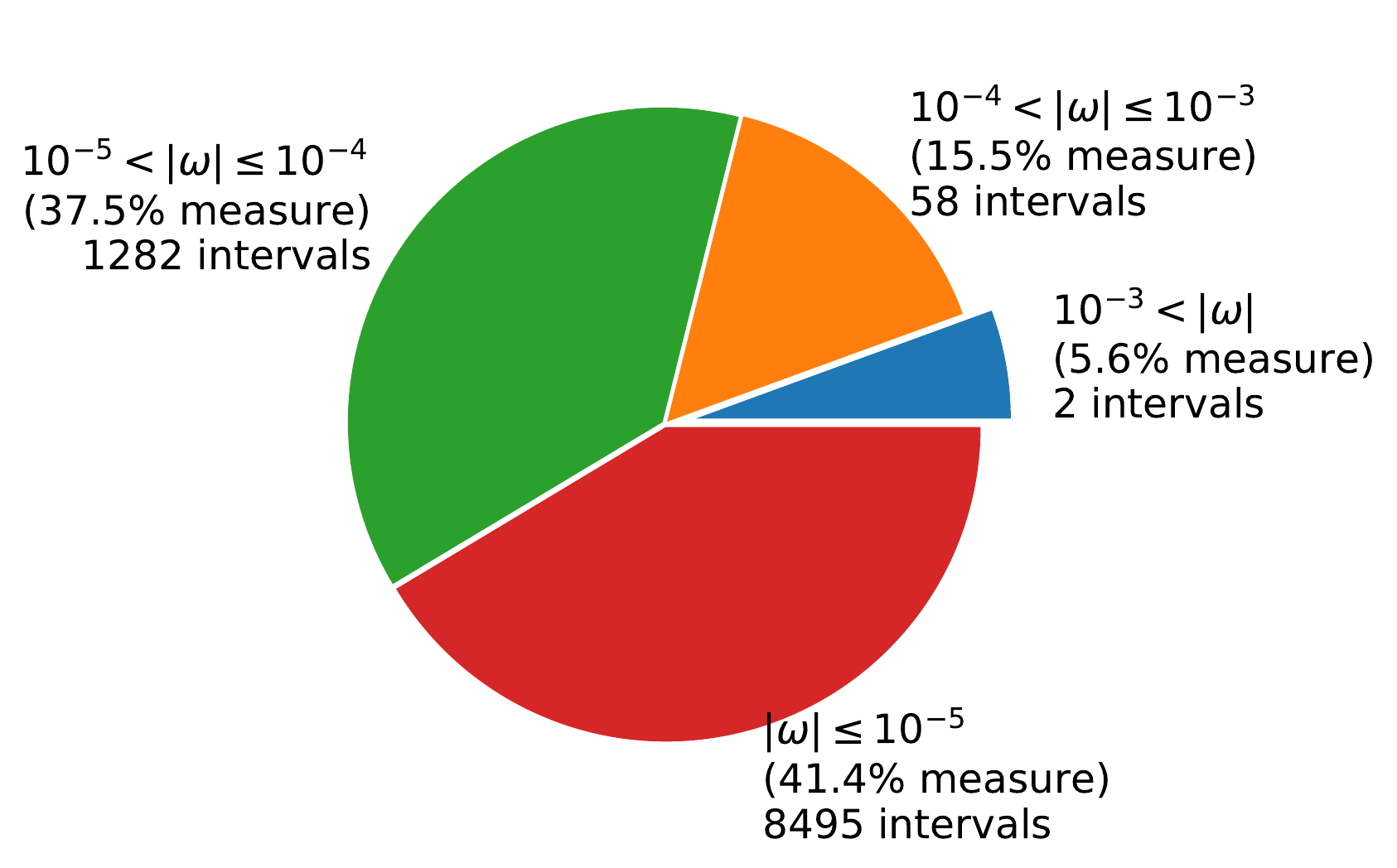}

\centerline{\hfill $\Omega_3 = [1.6,1.7]$, $\mu_3 \approx 92.4\%$ \hfill
\hfill $\Omega_4 = [1.7,1.8]$, $\mu_4 \approx 70.3\%$ \hfill}

\vskip 12pt

\noindent
\includegraphics[scale=.32]{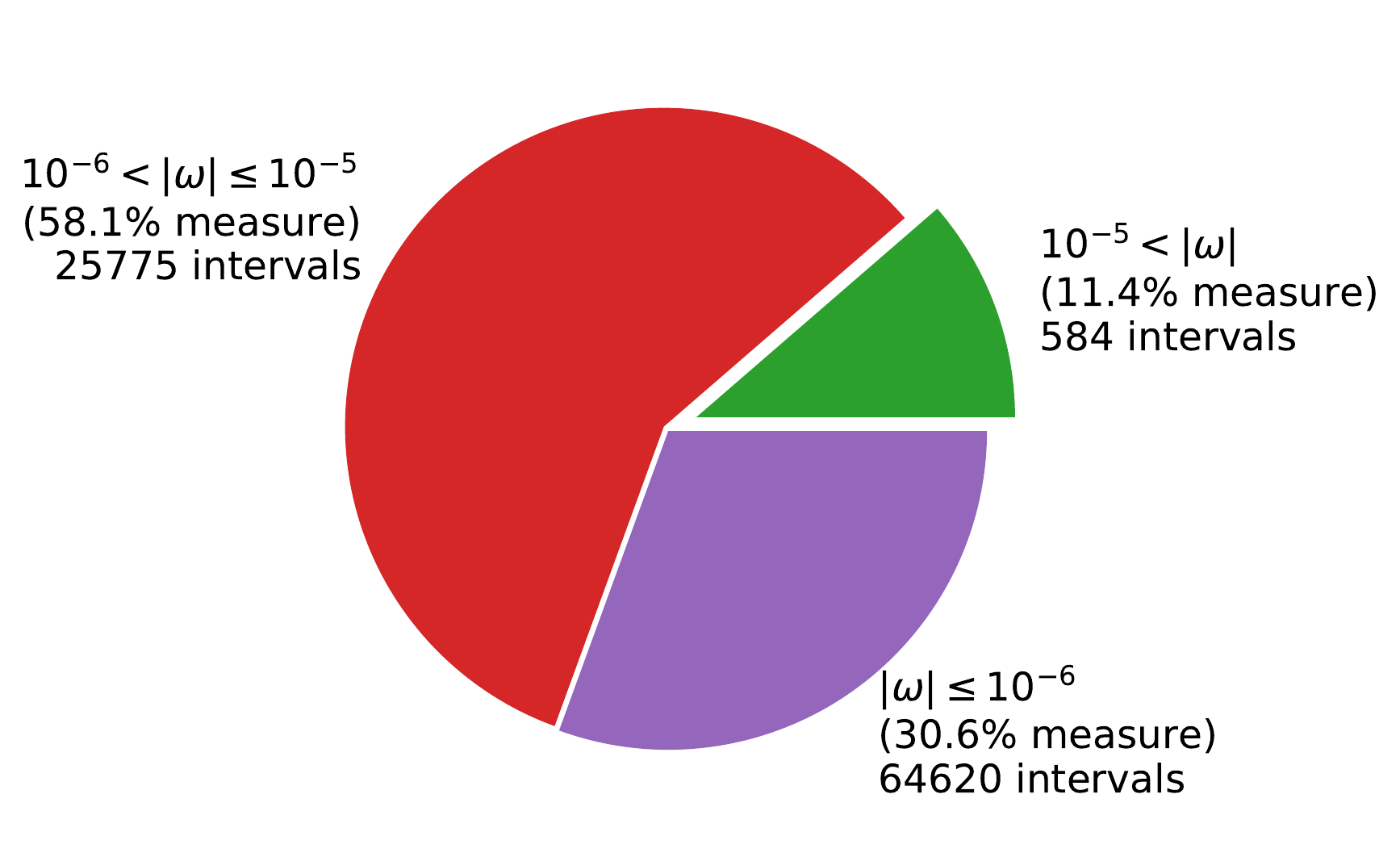}
\includegraphics[scale=.32]{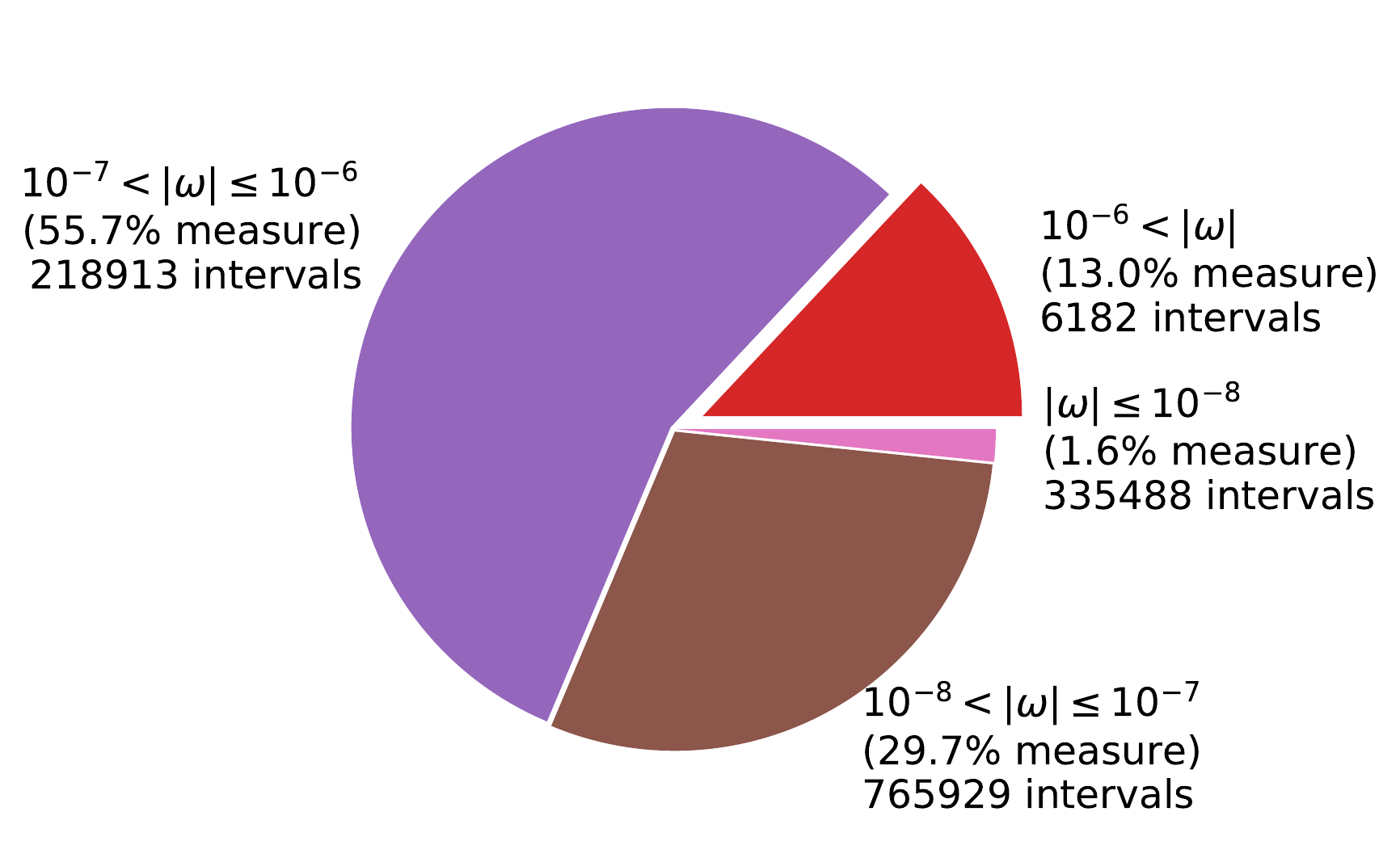}

\centerline{\hfill $\Omega_5 = [1.8,1.9]$, $\mu_5 \approx 97.5\%$ \hfill
\hfill $\Omega_6 = [1.9,2.0]$, $\mu_6 \approx 97.3\%$ \hfill}
\caption{\label{fig:sixPies}
Contribution of intervals of various widths to the total measure of $\cP^+$
in the smaller intervals $\Omega_1, \ldots, \Omega_6 \subset \Omega$.
The percentage $\mu_i$ of the measure of $\cP^+$ in each $\Omega_i$ is indicated.
Pie charts smaller than $2\%$ for small intervals were joined to the last one.
The first pie slice is slightly pulled out to indicate the contribution
of the largest intervals.
The colouring of width ranges is consistent among the six charts.}
\end{figure}

It turns out that the distribution of the intervals with the required escape time is very uneven across $\Omega$, so the result shown in \cite{GolKouLuzPil20} did not provide a full account of the situation. In order to show the proportions between the amounts of intervals of different sizes in different regions of the parameter space, now we subdivide $\Omega$ into six subsets $\Omega_1, \ldots, \Omega_6$ of the same length $0.1$, and we conduct the analogous computation for each subset separately. Figure~\ref{fig:sixPies} shows pie diagrams obtained for these six subintervals $\Omega_1, \ldots, \Omega_6$ of $\Omega$. Additionally, the total measure $\mu_i$ of the constructed intervals with the required escape time is given in terms of the percentage of the measure of each $\Omega_i$. It should not be surprising to see that this measure corresponds to the portion of parameters outside the periodic windows that can be spotted in the bifurcation diagram. Moreover, as we progress from $1.4$ to $2.0$ in $\Omega$, the contribution of gradually smaller intervals becomes more significant. The contrast is considerable and reflects the cost of the computation. For example, for $\Omega_1$, we obtain almost $400$ intervals and the computing time is below $25$ seconds. At the opposite end, for $\Omega_6$, we obtain over $1.3$ million intervals, and the computing time is well over $1$ hour.

A question arises on how much effort is actually necessary, in terms of analysing all the small intervals of parameters, in order to obtain large enough measure of the set of intervals in $\Omega$ with a satisfactory escape time. It is obvious that parameters in periodic windows must be excluded, and these comprise about $10\%$ of $\Omega$ ($0.06$ in measure), so in fact at most some $90\%$ of $\Omega$ ($0.54$ in measure) could be potentially covered by the intervals of interest. In order to answer this question, we conduct a series of complete computations with a few different values of $w$, ranging from $10^{-4}$ to $10^{-8}$, and with several different values of $N_0$, ranging from $15$ to $40$.

\begin{figure}[htbp]
\centerline{\includegraphics[width=10cm]{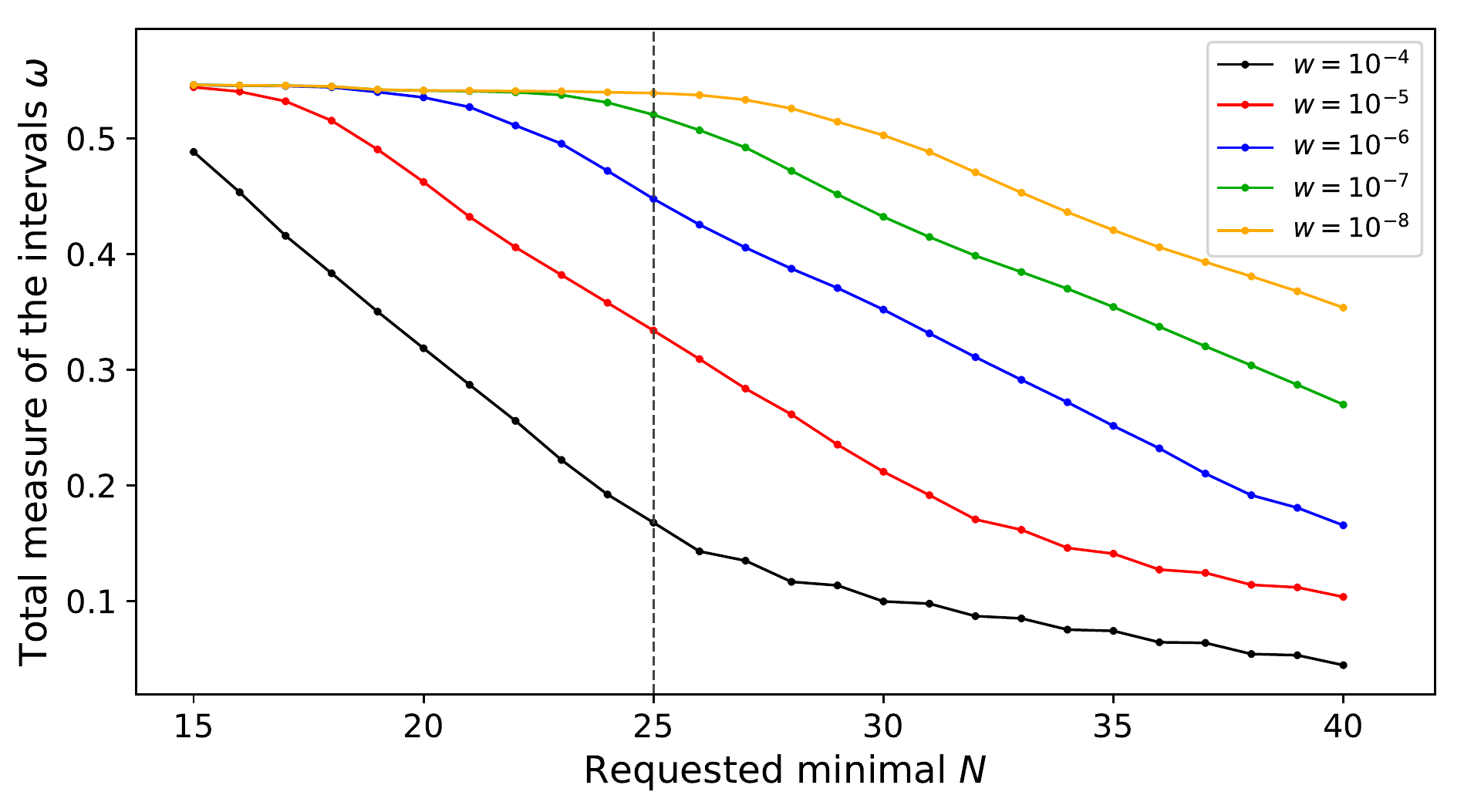}}
\caption{\label{fig:measureW} The total measure of subintervals of $\Omega$ that have an escape time $N(\omega) \geq N_0$, as a function of the requested time $N_0$. The computations were completed with different values of $w \in \{10^{-4}, \ldots, 10^{-8}\}$, and with $\delta = 10^{-3}$.}
\end{figure}

Figure \ref{fig:measureW} shows the results of these computations. For example, by looking at the points of intersection of the curves with the vertical dashed line, one can read from this figure that the total measure of parameter intervals with an escape time of at least $N_0 = 25$ is only about $0.16$ when the computations are conducted with $w = 10^{-4}$. This measure increases to some $0.33$ with $w$ decreased to $10^{-5}$, then it increases further to $0.44$ for $w = 10^{-6}$, then to $0.52$ for $w = 10^{-7}$, and eventually reaches some $0.539$ for $w = 10^{-8}$, which is nearly the maximum that could possibly be expected. The results are considerably better than shown in Figure~\ref{fig:noChNumMin_d3}, thanks to chopping off portions of the intervals that hit $\Delta$  and iterating the remaining portions or the intervals further.

It is a somewhat surprising observation that no matter how large $N$ is requested, the measure very close to $0.54$ can be apparently reached, provided the bound (defined by~$w$) on the allowed size of intervals is low enough. Unfortunately, this might result in very costly computations, because this measure must be filled up by a growing number of exponentially smaller intervals. For example, the computation time in the five discussed cases was ranging from some $5$ seconds for $w = 10^{-4}$ to about $25$ minutes for $w = 10^{-8}$.

\begin{figure}[htbp]
\centerline{\includegraphics[width=10cm]{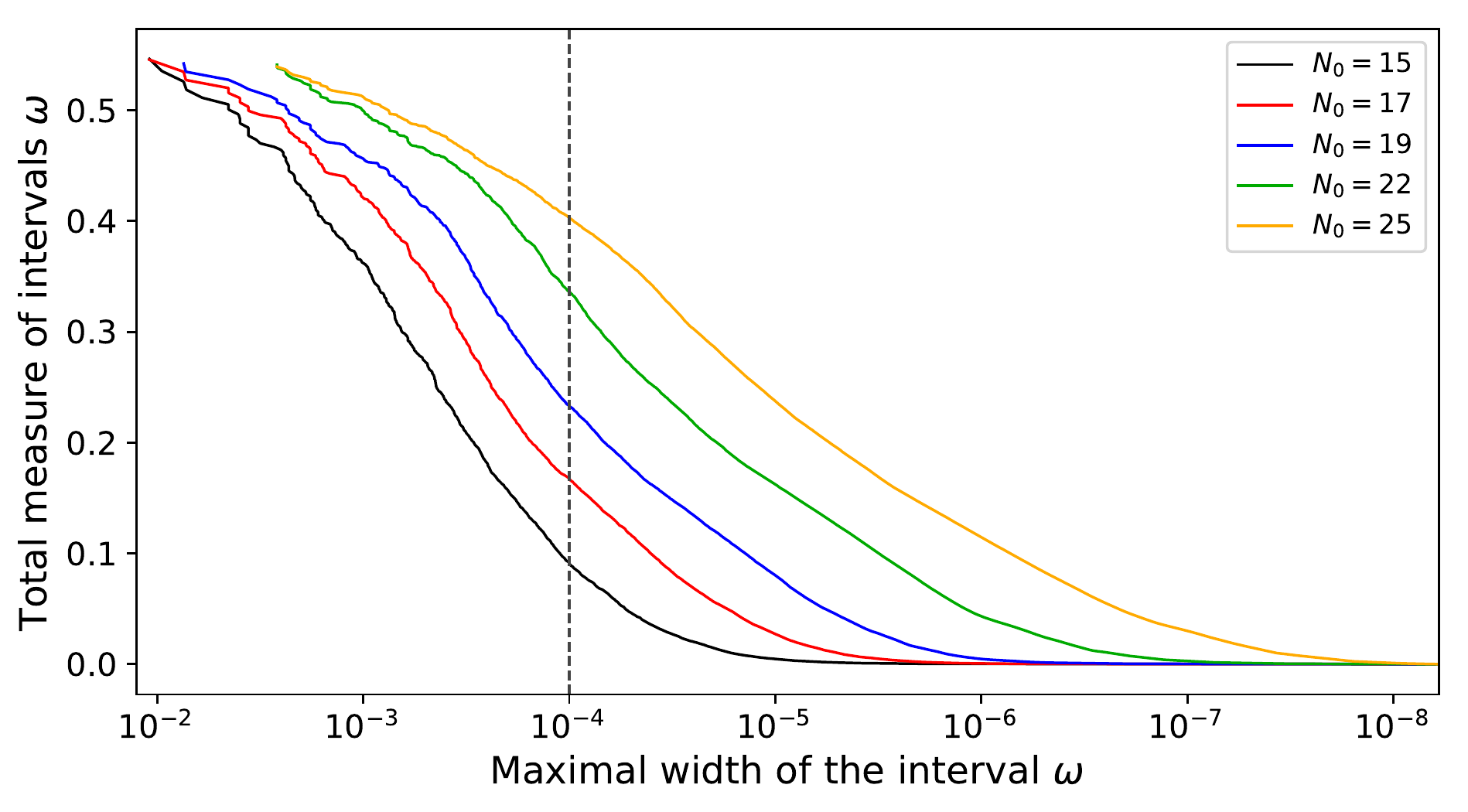}}
\caption{\label{fig:sizes8} The total measure of subintervals of $\Omega$ that have an escape time $N(\omega) \geq N_0$, shown for a few values of $N_0$, as a function of upper bound on the width of the subintervals. The computations were completed for $w = 10^{-8}$ and with $\delta = 10^{-3}$.}
\end{figure}

Figure \ref{fig:sizes8} illustrates the actual contribution of small intervals $\omega \subset \Omega$ to the total measure of intervals with an escape time at least $N_0$, obtained in the computations with $w = 10^{-8}$, shown for a few choices of $N_0$, selected on the basis of Figure~\ref{fig:measureW} in such a way that the total measure of these intervals is close to $0.54$. A point on each of the curves indicates the measure of intervals whose width is below the given threshold. For example, if one looks at the intersections of the curves with the dashed vertical line, one can see that the contribution of intervals narrower than $10^{-4}$ to the measure of $\cP^+$ in $\Omega$ is below $0.1$ for $N_0 = 15$, it is about $0.16$ for $N_0 = 17$, and is gradually increasing up to about $0.4$ for $N_0 = 25$. We did not plot these curves for higher values of $N_0$ because in those cases the measure $0.54$ was not reached due to the restrictive value of $w$, so the image would not be complete. It is worth to note the different starting points of the curves at the top of the plot; each of these points indicates the width of the widest interval with an escape time at least $N_0$. Moreover, the flattening bottom portions of the curves show that considering extremely small intervals provides gradually smaller gain in the total measure obtained.

\begin{figure}[htbp]
\centerline{\includegraphics[width=10cm]{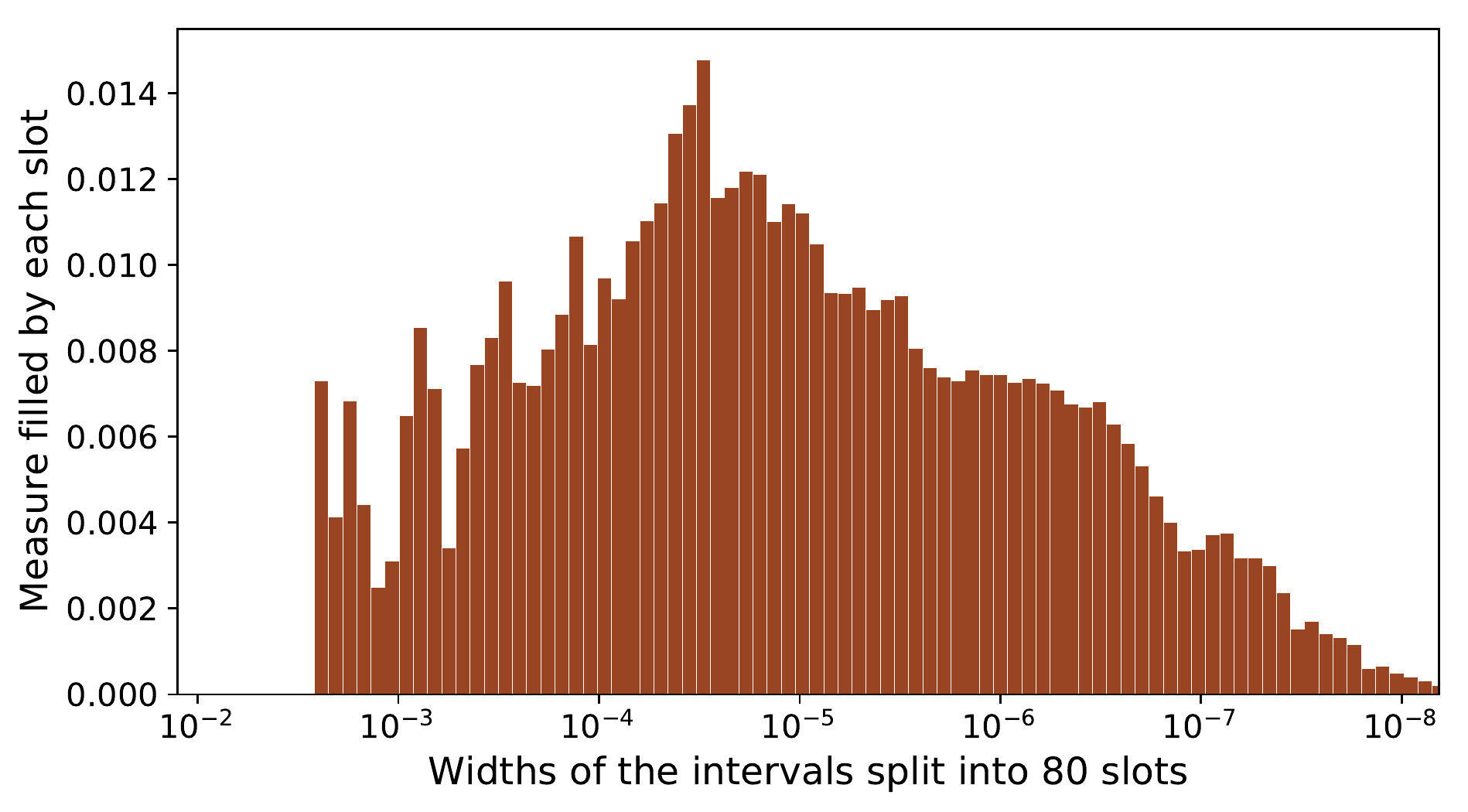}}
\caption{\label{fig:contrib825} The measure of subintervals of $\Omega$ that have an escape time $N(\omega) \geq N_0 = 25$ divided into slots depending on their widths. The computation was completed for $w = 10^{-8}$ and with $\delta = 10^{-3}$.}
\end{figure}

As a matter of fact, the contribution of intervals of various sizes to the total measure of intervals that have a satisfactory escape time can be illustrated in another way. In Figure~\ref{fig:contrib825}, all the intervals with an escape time of at least $25$ are gathered into 80 slots depending on their width. The slots are of equal size in the logarithmic scale, as shown in the diagram. One can clearly see the prevailing contribution of the intervals in the medium range of the widths. It appears in the picture that the ``tail'' of the slots was cut, and there might have been some intervals with sizes below $10^{-8}$ which were discarded due to the restrictive value of $w$. Indeed, in this computation, the total measure of the $1{,}222{,}230$ parameter intervals was proved to be at least $0.539302250926$, as opposed to the rigorous lower bound $0.539934844013$ on the measure obtained with $1{,}436{,}063$ parameter intervals constructed with $w = 10^{-10}$ in \cite{GolKouLuzPil20}; the difference is at least $6 \cdot 10^{-5}$.


\section{Choosing the number of bisection steps}
\label{sec:bisection}

The number $s$ of bisection steps affects the accuracy of chopping an interval $\omega$ for which $\omega_i \cap \Delta \neq \emptyset$ into smaller pieces that can be iterated further without hitting $\Delta$ at the $i$-th iterate. We refer to \cite[Section III~C and Algorithm 6.5]{GolKouLuzPil20} for the details of the method. In this section, we describe experimental study of the effect of this accuracy on the overall result of the computation.

\begin{figure}[htbp]
\centerline{\includegraphics[width=9cm]{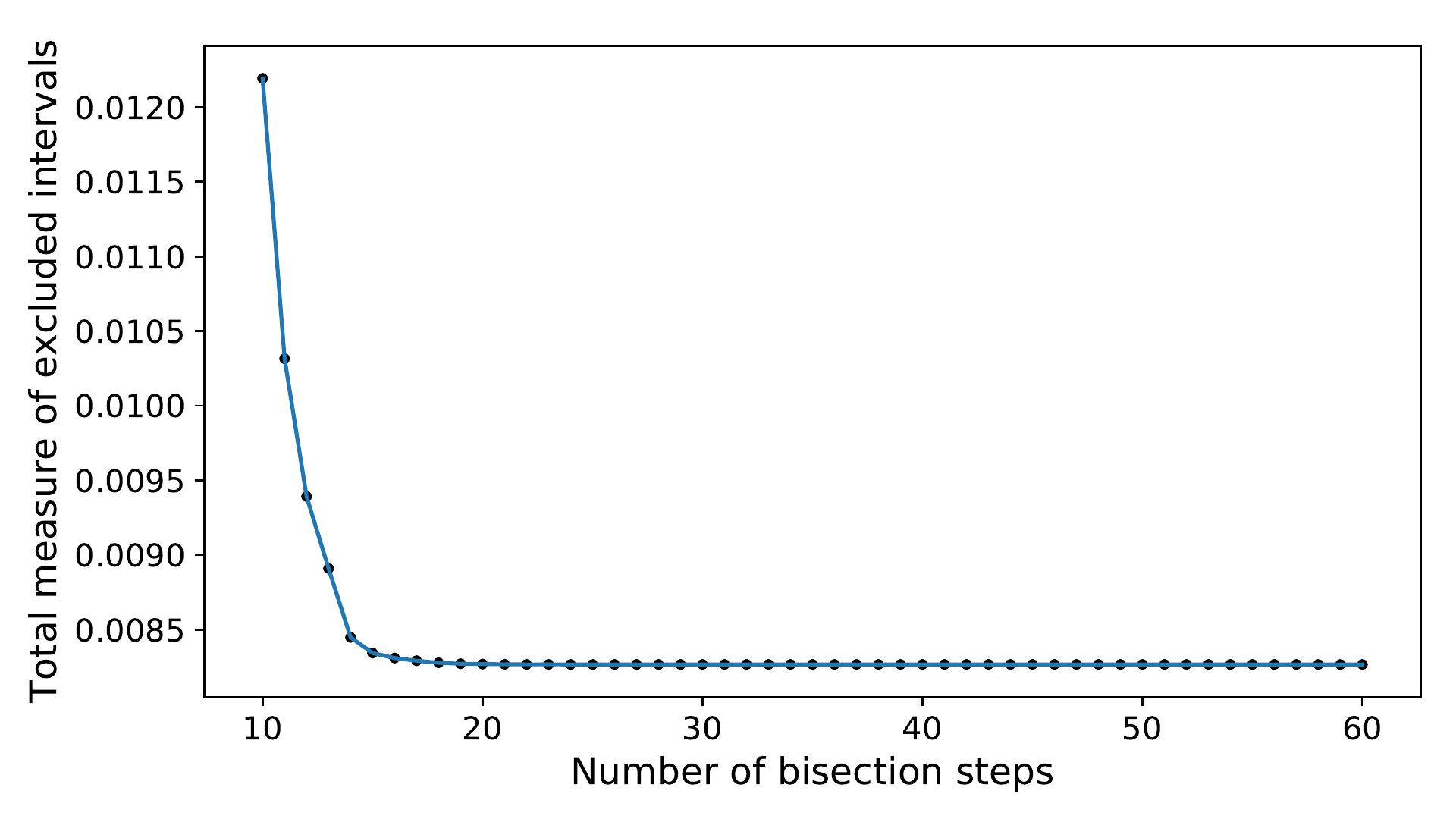}}
\caption{\label{fig:bisMeas10} The measure of parameters excluded due to the collision of some of their iterate with the critical neighbourhood $\Delta$, as a function of the number $s$ of bisection steps used to estimate the preimage of $\Delta$.}
\end{figure}

\begin{figure}[htbp]
\centerline{\includegraphics[width=9cm]{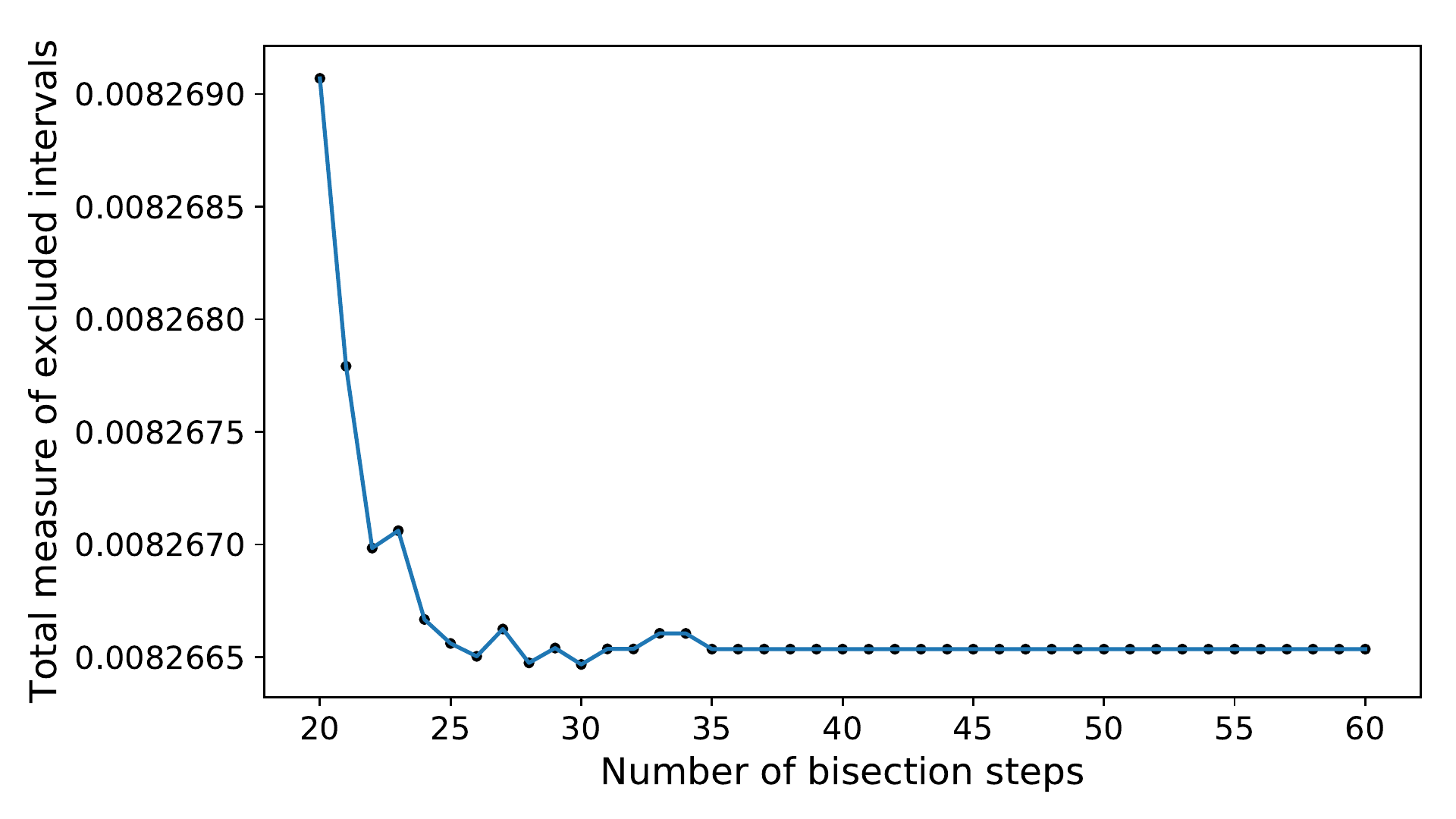}}
\caption{\label{fig:bisMeas20} The measure of parameters excluded due to the collision of some of their iterate with the critical neighbourhood $\Delta$, as a function of the number $s$ of bisection steps used to estimate the preimage of $\Delta$, shown for a narrower range of $s$ than in Figure~\ref{fig:bisMeas10}.}
\end{figure}

\begin{figure}[htbp]
\centerline{\includegraphics[width=9cm]{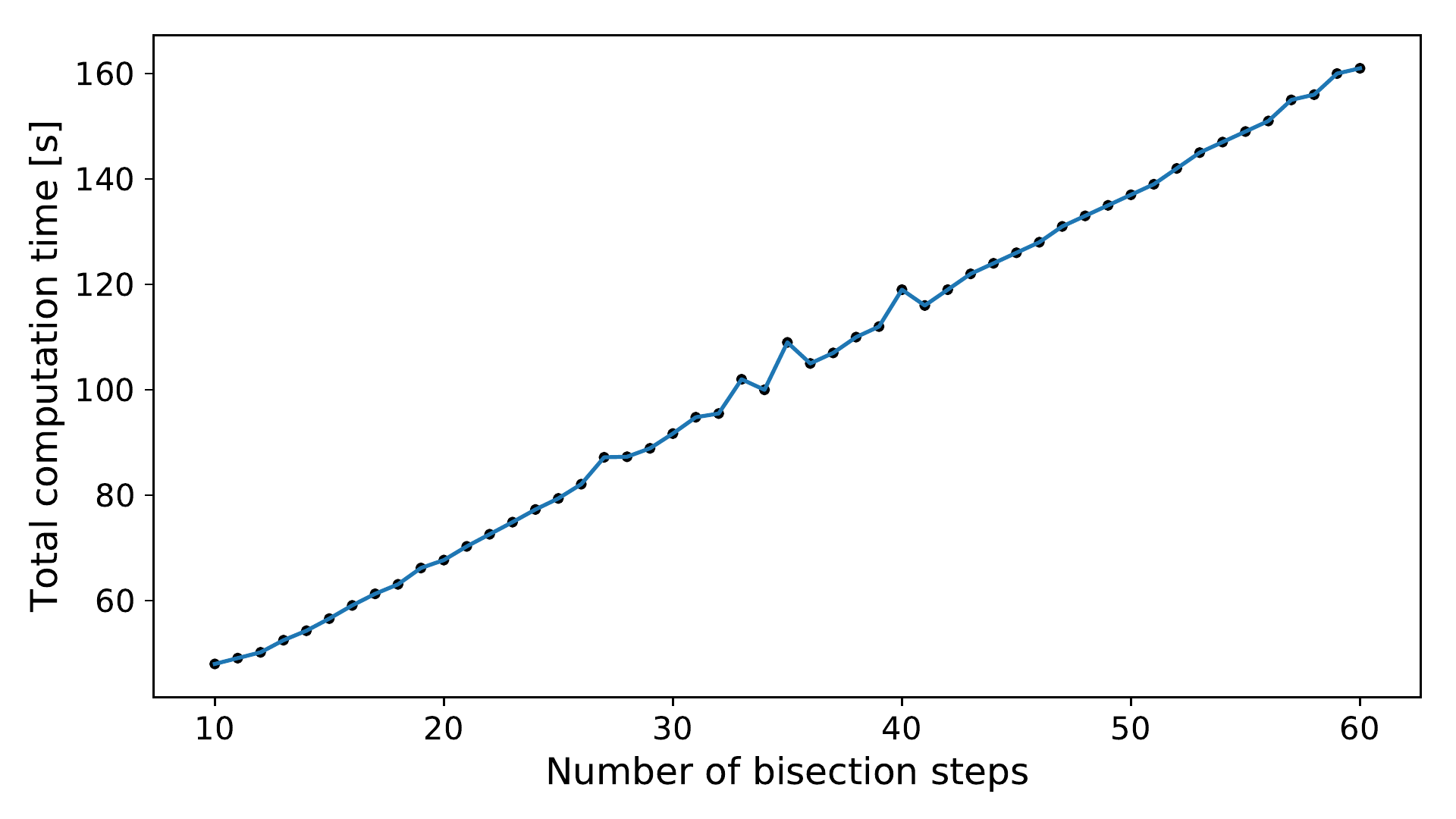}}
\caption{\label{fig:bisTime10} The computation time of iterating $100{,}000$ intervals, as a function of the number $s$ of bisection steps used to estimate the preimage of $\Delta$ at collision.}
\end{figure}

We iterate $100{,}000$ intervals starting with $\Omega = [1.4, 2]$, initially split uniformly into $u = 10$ subintervals, with $\delta = 10^{-3}$, $N_{\max} := N_0 := 100$, $w = 10^{-10}$, at the precision of $p = 1{,}000$ bits. We set the number $s = 10, 11, \ldots, 60$ of bisection steps at each attempt, and we compute the total measure of parameters excluded as an estimate for the preimage of $\Delta$ by $c_i$. Figure~\ref{fig:bisMeas10} shows rapid decrease with no further improvement after $s = 16$. However, a close-up shown in Figure~\ref{fig:bisMeas20} reveals some odd fluctuations that keep the result somewhat unstable until about $s = 36$. The time of computation increases linearly with the increase in the number $s$, as shown in Figure~\ref{fig:bisTime10}.

These calculations show that, on the one hand, increasing the number of bisection steps is not computationally expensive but, on the other hand, the gain is negligible after certain threshold. As a consequence, this number must be adjusted each time, depending on the specific calculations conducted.


\section*{Conclusion and final remarks}
\label{sec:final}

In the research reported on in the paper, we have found a comprehensive collection of intervals of parameters with certain rigorously proved dynamical properties. We plan to use these intervals in our further research aimed at the development of a full-featured computer-assisted method for computing a lower bound on the measure of stochastic parameters in the quadratic map, combined with other results already obtained \cite{DayKokLuzMisOkaPil08,GolKouLuzPil20,GolLuzPil16}, and those that still need to be completed. Thorough understanding of the dynamics that can be tracked with rigorous numerical methods is crucial for completing all the stages in the construction based in part on \cite{LuzTak06}. Moreover, the features observed in our computations may become a motivation for defining new notions that apply to dynamical models observed at finite scale and are motivated by the corresponding ``infinitesimal'' terms, see e.g.~\cite{LuzPil11}.

Last but not least, we would like to point out the unfortunate fact that it is a common practice in mathematical proofs, including computer-assisted rigorous numerical proofs, that the authors do not provide insight into how certain specific values of parameters were found or guessed, even though these values are often crucial for the success of the method applied. In contrast to this, our paper provides a systematic study of a wide range of adjustable settings; this study is aimed at finding those settings for which rigorous computation of non-zero measure of chaotic parameters has the best chances to succeed. We hope that our results and discussions shed light on how the various constraints and adjustments of the numerical method can be tweaked in order to achieve desired goals.



\end{document}